\documentclass{amsart}

\newtheorem{theorem}{Theorem}[section]
\newtheorem{lemma}[theorem]{Lemma}
\newtheorem{corollary}[theorem]{Corollary}

\theoremstyle{definition}
\newtheorem{definition}[theorem]{Definition}

\newtheorem{proposition}[theorem]{Proposition}

\theoremstyle{remark}
\newtheorem{remark}[theorem]{Remark}

\numberwithin{equation}{section}

\newcommand{\Real}{\operatorname{Re}}
\newcommand{\Imag}{\operatorname{Im}}
\newcommand{\HSC}{\operatorname{HSC}}
\newcommand{\R}{\mathbb R}
\newcommand{\Int}{\int_X}

\usepackage{geometry}
\begin{document}

\title[]{Holomorphic sectional curvature and spectral estimates on Hermitian surfaces}

\author[]{Jieming Yang}

\address{Department of mathematics, Wenzhou University, Zhejiang, China.}
\email{20230120@wzu.edu.cn}



\thanks{The author thanks Professor Jixiang Fu. The author also thanks ChatGPT for suggestions regarding the terms \eqref{eq:Sf-invariant} and \eqref{eq:Qstar}, the examples in Section 4.2, the literature search, the calculation checks, and the English writing. The author takes full responsibility for the manuscript.}

\date{}

\begin{abstract}
	On compact connected Gauduchon surfaces, we give a integral curvature criterion for a positive lower bound on the first positive eigenvalue of the Hodge Laplacian on functions. If the criterion holds, then for any prescribed complex line and any $B>0$, there is a perturbation in the same Aeppli class that retains at least half the original spectral lower bound and has holomorphic sectional curvature less than $-B$ on that line for every Gauduchon connection. We also extend the Wang–Yang eigenvalue estimate in complex dimension two to Gauduchon metrics under a directional lower bound for the holomorphic sectional curvature of a Gauduchon connection that includes an explicit torsion correction.
\end{abstract}

\maketitle

\section{Introduction}

Let $(M, g)$ be a closed Riemannian manifold of dimension $m\geq2$. The classical lower bounds for the first positive eigenvalue $\lambda_1$ of the Hodge Laplacian
\[
\Delta_{\mathrm d}=\mathrm d^{*}\mathrm d+\mathrm d\mathrm d^{*}
\]
begin with Ricci curvature. Lichnerowicz \cite{Lich} proved that if $\operatorname{Ric}\geq(m-1)\kappa\,g$ for some $\kappa>0$, then $\lambda_1\geq m\,\kappa$.

On a closed connected K\"ahler manifold $(X, \omega)$ of $\dim_{\mathbb{C}}X=n\geq2$, the corresponding estimate reads as follows (e.g., \cite{Futaki}): if $\operatorname{Ric}(\omega)\geq\kappa\,\omega$ for some $\kappa>0$, then $\lambda_1\geq2\,\kappa$.

Tsukamoto \cite{Tsukamoto} showed that a uniform positive lower bound for the holomorphic sectional curvature bounds the diameter of a complete K\"ahler manifold from above. Such a hypothesis need not give a positive Ricci lower bound, even on a compact surface.

Hitchin \cite{Hitchin} constructed K\"ahler metrics of positive holomorphic sectional curvature on the Hirzebruch surfaces; Yang and Zheng \cite{YangZheng} obtained such metrics in every K\"ahler class. Yet for $a\geq2$ the Hirzebruch surface $\mathbb F_a$ is not Fano and hence admits no K\"ahler metric of positive Ricci curvature. For these surfaces the classical Ricci estimate gives no positive lower bound, even though the holomorphic sectional curvature of Hitchin's metrics is uniformly positive.

Li and Wang \cite{LiWang} obtained an eigenvalue comparison on closed K\"ahler manifolds in terms of the diameter and separate lower bounds for holomorphic sectional curvature and orthogonal Ricci curvature.

Wang and Yang \cite{WY} derived a bound from holomorphic sectional curvature alone. Their integral identity retains the curvature of the complex line spanned by the gradient of an eigenfunction.

In complex dimension two their estimate reads: if the Chern holomorphic sectional curvature $\operatorname{HSC}^{0}\geq2$, then $\lambda_1\geq\frac{896}{225}$. Their product examples explain why the Fubini-Study eigenvalue is not the comparison value here: the metric of constant holomorphic sectional curvature $2$ on $\mathbb{CP}^2$ has first eigenvalue $6$, whereas on $\mathbb{CP}^1\times\mathbb{CP}^1$ there are product metrics with minimum holomorphic sectional curvature $2$ and first eigenvalue tending to $4$.

Let $(X, \omega)$ be a compact connected Hermitian surface. We study the corresponding question for the Gauduchon connections (\cite{Gauduchon97}):
\[
\nabla^t=(1-t)\,\nabla+t\,\nabla^{\mathrm B},\quad t\in\R,
\]
where $\nabla=\nabla^0$ is the Chern connection and $\nabla^{\mathrm B}$ is the Bismut connection (\cite{Bismut}). We fix the Chern connection as our background connection. If $\omega$ is not K\"{a}hler, these connections are all distinct. We refer to \cite{AOUV}, \cite{FZ}, \cite{LS}, \cite{WY1}, \cite{ZZ22}, \cite{ZZ23} and referencesin for this topic.

On a complex surface the Chern torsion is determined by the Gauduchon torsion (1,0)-form (\cite{Gauduchon}): $\tau=\Lambda_{\omega}\partial\omega$. The Gauduchon condition permits $\tau\ne0$ but makes it possible to cancel the relevant first-order term after integration.

Let $\HSC^t$ denote the holomorphic sectional curvature of the Gauduchon connection $\nabla^t$. For a $\varphi\in\mathcal{A}^{0}_{\mathbb{R}}(X)$, we write $\HSC^t\ge\varphi$ if $\HSC^t(v)\ge\varphi(x)$ for $x\in X$, $v\in(T_x^{1,0}X\backslash\{0\})$. Hereafter we may omit the volume form $\frac{\omega^{2}}{2!}$ for simplicity.

We consider how much of the Wang–Yang estimate survives for a (non-K\"{a}hler) Gauduchon metric, and whether curvature can still give a positive spectral bound when the pointwise hypothesis fails. Theorem \ref{thm:directional-spectral-gap} answers the first question. Theorem \ref{11062} and Proposition \ref{prop:nonclosed-local-negative} answer the second, even when the failure occurs on a prescribed complex line in a fixed Aeppli class.

\begin{theorem}\label{thm:directional-spectral-gap}
	Let $(X,\omega)$ be a compact connected Gauduchon surface. If there exists some $\kappa\in\mathbb R$, such that for every $v\in(T_x^{1,0}X\backslash\{0\})$,
	\begin{equation}\label{eq:directional-curvature-bound}
		\mathrm{HSC}^{0}(v)\geq\kappa+\frac1{25}\left(
		7\mid\!\tau\!\mid^2+8\frac{\mid\!\tau(v)\!\mid^2}{\mid\!v\!\mid^2}
		\right),
	\end{equation}
	then
	\[
	\lambda_1>\frac{448}{225}\,\kappa.
	\]
	More generally, let $f$ be a nonconstant smooth real solution to $\Delta_{\mathrm d}f=\psi(f)$, where $\psi$ is $C^{1}$ on a neighborhood of $f(X)$. Then
	\[
	\frac{\Int\psi'(f)\mid\!\nabla f\!\mid^2}{\Int\mid\!\nabla f\!\mid^2}
	>\frac{448}{225}\,\kappa.
	\]
\end{theorem}

\begin{remark}
	For any fixed $t\in\mathbb{R}$, the condition \eqref{eq:directional-curvature-bound} is equivalent to
	\[
	\HSC^{t}(v)
	\ge\kappa+\frac35\,\frac{\mid\!\tau(v)\!\mid^2}{\mid\!v\!\mid^2}+\left(\frac7{25}-t^2\right)\left(\mid\!\tau\!\mid^2-\frac{\mid\!\tau(v)\!\mid^2}{\mid\!v\!\mid^2}\right)
	\]
	for $x\in X$, $v\in(T^{1,0}_xX\setminus\{0\})$.
	
	Indeed, $\frac{\mid\tau(v)\mid^2}{\mid v\mid^2}$ is the squared norm of the restriction of $\tau$ to $\mathbb{C}v$, and the quantity in the last parentheses is the squared norm of its restriction to $\left(\mathbb{C}v\right)^\perp$. By \eqref{eq:rankone}, passing from $\HSC^0(v)$ to $\HSC^t(v)$ subtracts $t^2$ times this latter quantity.
	
	In particular, for either choice $t=\pm\frac{\sqrt{7}}{5}$, condition \eqref{eq:directional-curvature-bound} takes the simpler form
	\[
	\HSC^{t}(v)\ge\kappa+\frac35\,\frac{\mid\!\tau(v)\!\mid^2}{\mid\!v\!\mid^2}
	\]
	for $x\in X$, $v\in(T^{1,0}_xX\setminus\{0\})$. For these connections only the restriction of $\tau$ to the complex line $\mathbb Cv$ enters the correction.
	
	Thus the theorem extends the Wang-Yang eigenvalue estimate in complex dimension two to compact Gauduchon surfaces under the stated directional curvature bound. In the K\"ahler case, the torsion vanishes and the hypothesis reduces to $\HSC^0\ge\kappa$; the estimate retains the Wang-Yang coefficient and rules out equality.
\end{remark}

\begin{corollary} \label{cor:isotropic}
	Let $(X,\omega)$ be a compact connected Gauduchon surface. If there exist some $\kappa, t\in\R$, such that
	\begin{equation}\label{eq:main-pinching}
		\HSC^t\ge\kappa+\frac35\mid\!\tau\!\mid^2,
	\end{equation}
	then both conclusions of Theorem \ref{thm:directional-spectral-gap} hold. In particular, if $\HSC^t\ge\kappa$, then
	\[
	\lambda_1>\frac{448}{225}\left(\kappa-\frac35\mid\!\tau\!\mid_\infty^2\right),
	\]
	where $\displaystyle\mid\!\tau\!\mid_\infty^2=\max_{x\in X}\mid\!\tau\!\mid^2$.
\end{corollary}

\begin{corollary} \label{cor:closed-threshold}
	
	Fix $\kappa, t\in\mathbb R$. Let $(X, \omega)$ be a compact connected Gauduchon surface satisfying either \eqref{eq:directional-curvature-bound} or \eqref{eq:main-pinching}. Suppose that $f$ is a smooth real solution of $\Delta_\mathrm{d}f=\psi(f)$, where $\psi$ is $C^{1}$ on a neighborhood of $f(X)$. If
	\[
	\psi'(y)\le\frac{448}{225}\,\kappa
	\]
	for every $y\in f(X)$, then $f$ is a constant.
\end{corollary}

The estimate under the directional curvature bound and a second estimate under an integral criterion come from two choices of weight in the same formula. For Theorem \ref{thm:directional-spectral-gap}, we take $\left(\mid\!\nabla f\!\mid^2+\varepsilon\right)^{-1}$ and then let $\varepsilon\rightarrow0$. A drift term whose integral vanishes by the Gauduchon condition allows the transverse torsion term to be absorbed
into the same Hessian square as in the K\"ahler case. This accounts for the unchanged spectral coefficient. Equality would force the real Hessian at a maximum of a nonconstant first eigenfunction to have nonintegral rank. The bound is therefore strict for each metric; the
argument does not give a uniform increase in the coefficient. The semilinear estimate with weight $\mid\!\nabla f\!\mid^2$ follows by integration by parts and the variational characterization of $\lambda_1$.

With weight one, the K\"ahler case recovers the quartic eigenfunction identity of Wang-Yang in complex dimension two. On a Gauduchon surface, completing squares retains $\mid\nabla\mid\!\nabla f\!\mid^2\mid^2$ in the integral.

Define
\begin{align}
	\mathcal H_\omega(x)
	& =\min_{v\in T_x^{1,0}X,\, \mid v\mid=1}\HSC^{0}(v)-\frac12\mid\!\tau\!\mid^2(x), \label{eq:effective-quartic-potential}\\
	\sigma_\omega=\lambda_{\min}\!\left(\frac14\Delta_{\mathrm d}+\mathcal H_\omega\right)
	& = \inf_{u\in(W^{1,2}(X;\mathbb R))\setminus\{0\}}
	\frac{\int_X\left(\tfrac12 \mid\!\nabla u\!\mid^2+\mathcal H_\omega\,u^2\right)}
	{\int_Xu^2}.\label{eq:curvature-form-bottom}
\end{align}
where $\HSC^{0}(v)$ is given in Definition \ref{11161}, and $W^{1,2}(X;\mathbb R)$ denotes the Sobolev space of real-valued $L^2$ functions with weak differential in $L^2$. Theorem \ref{11062} gives $\lambda_1>\sigma_\omega$ on every compact connected Gauduchon surface. Moreover, every nonconstant smooth real solution of $\Delta_{\mathrm d}f=\psi(f)$, with $\psi\in C^1$ near $f(X)$, satisfies
\[
\frac{\int_X \psi'(f)\mid\!\nabla f\!\mid^4}{\int_X\mid\!\nabla f\!\mid^4}>\sigma_\omega.
\]

Proposition \ref{prop:nonclosed-local-negative} shows that this curvature lower bound can survive a local loss of positivity. Suppose that $\sigma_{\omega_0}>0$, and fix a point $p$, a complex line $L\subset T^{1,0}_pX$, and a number $B>0$. In the Aeppli class of $\omega_0$ there are Gauduchon metrics $\omega_\varepsilon$ converging to $\omega_0$ in $C^1$, with the same exterior derivative and total volume, for which
\[
\HSC^{t}_{\omega_\varepsilon}(L)<-B, \quad
\lambda_1(\omega_\varepsilon)>\sigma_{\omega_\varepsilon}
\ge\tfrac12\sigma_{\omega_0}>0,
\]
where $t\in\R$. Thus the directional hypothesis fails for every positive $\kappa$, while the integral criterion still gives a positive lower bound on the same metric. What persists under this perturbation is the lower bound defined by curvature, not merely the positivity of the first eigenvalue.

The distinction between the two curvature criteria concerns the metric. If a compact complex surface admits a Gauduchon metric satisfying either of the two positive criteria, then it is rational (Corollary \ref{cor:gs-positive-criteria-rationality}). The converse uses Zhang's theorem that every rational surface admits a K\"ahler metric of positive holomorphic sectional curvature \cite[Theorem 1.4]{ShiZhang}. Perturbing such a metric gives a (non-K\"{a}hler) Gauduchon metric for which both bounds are positive. In particular, the given metric need not be K\"ahler, although its underlying complex surface must admit a K\"ahler metric.

Two comparisons help place these estimates. Zhang \cite{Zhang} uses a lower bound for Bismut holomorphic sectional curvature on balanced manifolds together with a torsion commutator condition along a first eigendirection. On a surface, balancedness already implies K\"{a}hlerness. The Gauduchon hypothesis used here permits non-K\"{a}hler metrics and requires no separate condition along an eigendirection. Pure and mixed complex Hessians also enter the K\"ahler Sobolev inequalities studied by Baudoin and Munteanu \cite{BM}, Chakraborty \cite{Chak}, and Chakraborty, Wang and Wu \cite{CWW}. In those works the curvature hypothesis is a Ricci bound, while the curvature in our identities is the holomorphic sectional curvature in the gradient direction.

The paper is organized as follows. Section \ref{11072} fixes the curvature conventions and records the Gauduchon integration formulas. Section \ref{11073} proves the integral estimate and then the estimate under the directional curvature bound. Section \ref{sec:nonclosed-curvature-examples} separates the two curvature criteria by perturbing the metric in a fixed Aeppli class. A cyclic-cover argument yields the rationality conclusion and, when applied to a homogeneous Hopf metric, shows that a positive lower bound for $\HSC^{0}$ alone does not give a uniform spectral gap.

\section{Preliminaries}\label{11072}
Let $X$ be a complex manifold with $\dim_{\mathbb{C}}X=n\geq2$. Let
\[
\omega=\sqrt{-1}g_{i\bar{j}}\mathrm{d}z^{i}\wedge\mathrm{d}\bar{z}^{j}
\]
be a Hermitian metric on $X$.

Denote the complex and the Hodge Laplacian of $\omega$ by
\[
\Delta=\sqrt{-1}\Lambda_{\omega}\partial\bar{\partial}, \quad
\Delta_{\mathrm{d}}=\mathrm{d}\mathrm{d}^{\ast}+\mathrm{d}^{\ast}\mathrm{d},
\]
respectively.

Denote
\[
\Gamma^{p}_{ik}=g^{p\bar{j}}\partial_{k}g_{i\bar{j}}.
\]
The Chern connection is the unique Hermitian connection on $T^{1,0}X$ whose $(0,1)$-part is $\bar\partial$. Locally, we have
\[
\nabla_{k}\left(\frac{\partial}{\partial z^{i}}\right)
=\Gamma^{p}_{ik}\frac{\partial}{\partial z^{p}}, \quad
\nabla_{\bar{l}}\left(\frac{\partial}{\partial z^{i}}\right)=0,
\]
where $\nabla_{k}=\nabla_{\frac{\partial}{\partial z^{k}}}$ and $\nabla_{\bar{l}}=\nabla_{\frac{\partial}{\partial\bar{z}^{l}}}$ denote the covariant derivatives of $\nabla$ along the vector fields $\frac{\partial}{\partial z^{k}}$ and $\frac{\partial}{\partial\bar{z}^{l}}$.

The Chern connection $\nabla$ serves as the background among the Gauduchon connections $\nabla^{t}$ in our later calculations. 

For every $\theta=\theta_{\bar{l}}\mathrm{d}\bar{z}^{l}\in\mathcal{A}^{0,1}(X)$, we have
\[
\partial_{k}\theta_{\bar{l}}=\nabla_{k}\theta_{\bar{l}},\quad \partial_{\bar{l}}\bar{\theta}_{\bar{k}}=\nabla_{\bar{l}}\bar{\theta}_{\bar{k}}.
\]
These terms may be used interchangeably.

The torsion tensor of the Chern connection $\nabla$ is given by
\[
T=T_{ki\bar{j}}\mathrm{d}z^{k}\otimes\mathrm{d}z^{i}\otimes\mathrm{d}\bar{z}^{j},
\]
where
\begin{equation}\label{11033}
	T_{ki\bar{j}}=\partial_{k}g_{i\bar{j}}-\partial_{i}g_{k\bar{j}}.
\end{equation}
Denote
\[
T^{p}_{ki}=g^{p\bar{j}}T_{ki\bar{j}}=\Gamma^{p}_{ik}-\Gamma^{p}_{ki}
\]
and
\[
T_{k}=T^{i}_{ki}=g^{i\bar{j}}T_{ki\bar{j}}.
\]
Then the torsion (1,0)-form is given by
\[
\tau=\Lambda_{\omega}\partial\omega=T_{k}\mathrm{d}z^{k}.
\]

Since each Hermitian metric $\omega$ satisfies
\begin{equation}\label{06041}
	\partial\omega^{n-1}=\tau\wedge\omega^{n-1},
\end{equation}
$\omega$ is Gauduchon if and only if (e.g., section 7.1 in \cite{LT})
\begin{equation}\label{10111}
	0=\mid\!\tau\!\mid^{2}+g^{k\bar{l}}\nabla_{\bar{l}}T_{k}.
\end{equation}
where
\[
\mid\!\tau\!\mid^{2}=g^{k\bar{l}}T_{k}\bar{T}_{l}, \quad g^{k\bar{l}}\nabla_{\bar{l}}T_{k}=g^{k\bar{l}}\nabla_{k}\bar{T}_{l}.
\]
In particular, for $n=2$, we obtain from \eqref{11033} and \eqref{06041} that
\begin{equation}\label{11034}
	T_{ki\bar{j}}=T_{k}g_{i\bar{j}}-T_{i}g_{k\bar{j}}.
\end{equation}

Integration by parts will be used frequently. Let $(X, \omega)$ be a compact Hermitian manifold. For $\theta=\theta_{\bar{l}}\mathrm{d}\bar{z}^{l}\in\mathcal{A}^{0,1}(X)$, we obtain from \eqref{06041} that
\[
\begin{split}
	0
	& = n\int_{X}\sqrt{-1}\mathrm{d}\left(\left(\theta-\bar{\theta}\right)\wedge\omega^{n-1}\right)\\
	& = n\int_{X}\sqrt{-1}\left(\partial\theta-\bar{\partial}\bar{\theta}+\tau\wedge\theta+\bar{\theta}\wedge\bar{\tau}\right)\wedge\omega^{n-1}\\
	& = \int_{X}2\operatorname{Re}\left(g^{k\bar{l}}\left(\nabla_{k}\theta_{\bar{l}}+T_{k}\theta_{\bar{l}}\right)\right)\omega^{n}.
\end{split}
\]

The curvature $R_{\omega}$ of the Chern connection $\nabla$ of a Hermitian metric $\omega$ is
\[
R_{\omega}=\nabla\circ\nabla\in\mathcal{A}^{1,1}\left(\mathrm{End}\left(T^{1,0}X\right)\right).
\]
Locally, we calculate
\[
R_{\omega}=\bar{\partial}\left(\partial gg^{-1}\right)^{T}
=R^{p}_{ik\bar{l}}\,\mathrm{d}z^{k}\wedge\mathrm{d}\bar{z}^{l}\otimes\mathrm{d}z^{i}\otimes\frac{\partial}{\partial z^{p}},
\]
where (see page 12 in \cite{Kobayashi})
\begin{equation}\label{09221}
	R^{p}_{ik\bar{l}}=g^{p\bar{j}}(-\partial_{\bar{l}}\partial_{k}g_{i\bar{j}}+g^{m\bar{n}}\partial_{\bar{l}}g_{m\bar{j}}\partial_{k}g_{i\bar{n}}).
\end{equation}
Denote
\[
R_{i\bar{j}k\bar{l}}=g_{p\bar{j}}R^{p}_{ik\bar{l}}.
\]
Then we have
\begin{equation}\label{10131}
	R_{i\bar{j}k\bar{l}}=\bar{R}_{j\bar{i}l\bar{k}}.
\end{equation}
The following identities are well-known (e.g., (2.6) in \cite{SW}):
\begin{equation}\label{11022}
	R_{i\bar{j}k\bar{l}}-R_{k\bar{j}i\bar{l}}
	=\nabla_{\bar{l}}T_{ik\bar{j}}, \quad
	R_{i\bar{j}k\bar{l}}-R_{k\bar{l}i\bar{j}}
	=\nabla_{\bar{j}}T_{ik\bar{l}}+\nabla_{k}\bar{T}_{jl\bar{i}}.
\end{equation}

The first and second Ricci curvatures of the Chern connection $\nabla$ are given by (e.g., \cite{Kobayashi})
\[
R_{i\bar{j}}=g^{k\bar{l}}R_{k\bar{l}i\bar{j}}, \quad 
K_{i\bar{j}}=g^{k\bar{l}}R_{i\bar{j}k\bar{l}}.
\]
To each of these, we can associate a real (1, 1)-form
\[
\rho_{\nabla}=\sqrt{-1}R_{i\bar{j}}\mathrm{d}z^{i}\wedge\mathrm{d}\bar{z}^{j}, \quad
\kappa_{\nabla}=\sqrt{-1}K_{i\bar{j}}\mathrm{d}z^{i}\wedge\mathrm{d}\bar{z}^{j}.
\]
There are two scalar curvatures:
\[
s=g^{i\bar{j}}g^{k\bar{l}}R_{i\bar{j}k\bar{l}}, \quad
\hat{s}=g^{i\bar{l}}g^{k\bar{j}}R_{i\bar{j}k\bar{l}}.
\]

The Bismut connection $\nabla^{B}$ of $\omega$ is defined by (\cite{Bismut})
\[
\nabla^{B}_{k}\left(\frac{\partial}{\partial z^{i}}\right)
=\Gamma^{p}_{ki}\frac{\partial}{\partial z^{p}},\quad
\nabla^{B}_{\bar{l}}\left(\frac{\partial}{\partial z^{i}}\right)
=g^{p\bar{q}}\bar{T}_{lq\bar{i}}\frac{\partial}{\partial z^{p}}.
\]
Then for $t\in\mathbb{R}$, the Gauduchon connection $\nabla^{t}$ is defined by (e.g., \cite{WY1})
\[
\nabla^{t}_{k}\left(\frac{\partial}{\partial z^{i}}\right)
=\left((1-t)\Gamma^{p}_{ik}+t\,\Gamma^{p}_{ki}\right)\frac{\partial}{\partial z^{p}},\quad
\nabla^{t}_{\bar{l}}\left(\frac{\partial}{\partial z^{i}}\right)
=t\,g^{p\bar{q}}\bar{T}_{lq\bar{i}}\frac{\partial}{\partial z^{p}}.
\]

The curvature $F^{t}$ of $\nabla^{t}$ is
\[
F^{t}= \left(\frac{1}{2}F^{p}_{imk}\mathrm{d}z^{m}\wedge\mathrm{d}z^{k}+\frac{1}{2}F^{p}_{i\bar{n}\bar{l}}\mathrm{d}\bar{z}^{n}\wedge\mathrm{d}\bar{z}^{l}
+F^{p}_{ik\bar{l}}\mathrm{d}z^{k}\wedge\mathrm{d}\bar{z}^{l}\right)\otimes\frac{\partial}{\partial z^{^{p}}}\otimes\mathrm{d}z^{i},
\]
where (e.g., \cite{FZ}, \cite{WY1})
\[
F^{p}_{ik\bar{l}}= R^{p}_{ik\bar{l}}+t\,g^{p\bar{j}}\left(\nabla_{\bar{l}}T_{ki\bar{j}}+\nabla_{k}\bar{T}_{lj\bar{i}}\right) +t^{2}\,g^{p\bar{j}}g^{m\bar{n}}\left(T_{ki\bar{n}}\bar{T}_{lj\bar{m}}-T_{mk\bar{j}}\bar{T}_{nl\bar{i}}\right).
\]
Denote
\[
F^{t}_{i\bar{j}k\bar{l}}=g_{p\bar{j}}F^{p}_{ik\bar{l}}.
\]
Then we have
\begin{equation}\label{eq:gaud-curvature-general}
	F^{t}_{i\bar{j}k\bar{l}}
	=R_{i\bar{j}k\bar{l}}+t\,\left(\nabla_{\bar{l}}T_{ki\bar{j}}+\nabla_{k}\bar{T}_{lj\bar{i}}\right) +t^{2}\,g^{p\bar{q}}\left(T_{ki\bar{q}}\bar{T}_{lj\bar{p}}-T_{pk\bar{j}}\bar{T}_{ql\bar{i}}\right).
\end{equation}
By \eqref{10131}, we have
\[
F^{t}_{i\bar{j}k\bar{l}}=\bar{F}^{t}_{j\bar{i}l\bar{k}}.
\]

\begin{definition}\label{11161}
	Let $(X, \omega)$ be a Hermitian manifold. The holomorphic sectional curvature of the Gauduchon connection along the tangent vector $v\in(T_x^{1,0}X\backslash\{0\})$ is
	\[
	\HSC^t(v)=\frac{F^t\left(v,\bar v,v,\bar v\right)}{\mid\!v\!\mid^4}.
	\]
\end{definition}

In particular, on a complex surface, for $f\in\mathcal{A}^{0}_{\R}(X)$ and $t\in\R$, we obtain from \eqref{11034}, \eqref{eq:gaud-curvature-general} and the above definition that
\begin{equation}\label{eq:rankone}
	F^t_{i\bar j k\bar l}f^i\bar f^{j}f^k\bar f^{l}
	=R_{i\bar j k\bar l}f^i\bar f^{j}f^k\bar f^{l}
	+t^2\left(\mid\!\nabla f\!\mid^2\mid\!f^iT_i\!\mid^2-\mid\!\tau\!\mid^2\mid\!\nabla f\!\mid^4\right).
\end{equation}

\begin{definition}\label{10112}
	Let $(X, \omega)$ be a Hermitian manifold. For $\ell\in\{1,2\}$, the $\ell$-Gauduchon-Ricci curvature form of $\nabla^{t}$ is
	\[
	F^{\ell}_{t}=\sqrt{-1}F^{\ell}_{i\bar{j}}\mathrm{d}z^{i}\wedge\mathrm{d}\bar{z}^{j},
	\]
	where
	\[
	F^{1}_{i\bar{j}}=g^{k\bar{l}}F^{t}_{k\bar{l}i\bar{j}}, \quad
	F^{2}_{i\bar{j}}=g^{k\bar{l}}F^{t}_{i\bar{j}k\bar{l}}.
	\]
\end{definition}

For simplicity, we will use the following notation. For $f\in\mathcal{A}^{0}_{\R}(X)$, set
\[
\begin{gathered}
	f_i=\partial_i f,\quad f_{\bar j}=\partial_{\bar j}f,\quad
	f^i=g^{i\bar j}f_{\bar j},\quad
	\bar f^{j}=g^{i\bar j}f_i, \quad \mid\!\nabla f\!\mid^2=g^{i\bar j}f_i f_{\bar j}=\frac{1}{2}\mid\!\mathrm{d}f\!\mid^2,\\
	f_{i\bar j}=\partial_i\partial_{\bar j}f=\nabla_{i}f_{\bar j}=\nabla_{\bar j}f_{i},\quad
	\mid\!f_{i\bar{l}}\!\mid^2=g^{i\bar j}g^{k\bar l}f_{i\bar{l}}f_{k\bar{j}}, \quad
	\mid\!\nabla_{k}f_{i}\!\mid^2
	=g^{i\bar j}g^{k\bar l}\nabla_{k}f_{i}\nabla_{\bar{l}}f_{\bar{j}},\\
	\mid\!f^kf_{k\bar{j}}\!\mid^2=g^{i\bar j}f^k\bar{f}^l f_{k\bar{j}}f_{i\bar{l}}, \quad
	\mid\!f^k\nabla_k f_i\!\mid^2=g^{i\bar j}f^k\bar f^{l}
	\nabla_k f_i \nabla_{\bar l}f_{\bar j}.
\end{gathered}
\]
On complex surface, by \eqref{11034}, we have
\begin{equation}\label{eq:bar-hessian-commute}
	\nabla_k f_i-\nabla_i f_k=T_i f_k-T_k f_i, \quad
	\nabla_{\bar l}f_{i\bar j}-\nabla_{\bar j}f_{i\bar l}
	=f_{i\bar l}\bar{T}_{j}-f_{i\bar j}\bar{T}_{l}.
\end{equation}

Denote
\[
\Phi(f)=2\operatorname{Re}\left(g^{k\bar{l}}T_{k}f_{\bar{l}}\right).
\]
Then we have (e.g., Lemma 7.1.1 in \cite{LT})
\[
\Delta_{\mathrm{d}}f=-2\Delta f-\Phi(f).
\]
\begin{lemma}\label{11041}
	Let $(X, \omega)$ be a compact Hermitian surface. For every $f\in\mathcal{A}^{0}_{\R}(X)$, the following identity holds
	\[
	0=\Int\left(\mid\!f_{i\bar j}\!\mid^2-\left(\Delta f\right)^2\right).
	\]
	Moreover, if $\omega$ is a Gauduchon metric, then for every $h\in\mathcal{A}^{0}_{\R}(X)$, the following identity holds
	\[
	\Int\Phi\left(h\right)=0.
	\]
\end{lemma}
\begin{proof}
	The first identity follows from (since $\dim_{\mathbb{C}}X=2$)
	\[
	0=\Int\sqrt{-1}\partial\bar{\partial}f\wedge\sqrt{-1}\partial\bar{\partial}f.
	\]
	While the second one follows from \eqref{10111}.
\end{proof}

Only the second identity requires the Gauduchon condition. We will extend the integrands by zero off
\[
X_+=\left\{\mid\!\nabla f\!\mid^2>0\right\}.
\]
The first identity may also be integrated over $X_+$, since both Chern Hessians of $f$ vanish almost everywhere on its critical set. Indeed, at every density point of the critical set, the differential of each real coordinate derivative of $f$ vanishes; the connection terms vanish there as well because $\mathrm{d}f=0$.

\section{Proofs}\label{11073}

In this section, we will prove our main results. The proofs rely on careful computations of the torsion, as well as two choices of weight in Lemma \ref{lem:commutator}.
\begin{lemma}\label{lem:commutator}
	Let $(X, \omega)$ be a compact Hermitian surface. Let $I$ be an open interval containing the range of $\mid\!\nabla f\!\mid^{2}$. For every $f\in\mathcal{A}^{0}_{\mathbb{R}}(X)$, $\chi=\chi(\mid\!\nabla f\!\mid^{2})\in C^{1}\left(I, \mathbb{R}\right)$, the following identity holds
	\begin{equation}\label{eq:generalweight}
		\begin{split}
			\Int\chi\,R_{i\bar j k\bar l}f^i\bar f^{j}f^k\bar f^{l}
			& = \Int\chi\,\left(\mid\!f^{k}f_{k\bar j}\!\mid^2-\mid\!f^k\nabla_k f_i\!\mid^2\right)\\
			& \quad +\Int\chi\,\operatorname{Re}\left(g^{i\bar j}f^k\nabla_kf_i\left(f^p f_{p\bar j}-\bar f^{l}\nabla_{\bar j}f_{\bar l}\right)\right)\\
			& \quad +\frac12 \Int\chi\left(-f^i\bar f^{j}f_{i\bar j}
			+\operatorname{Re}\left(f^if^k\nabla_i f_k\right)\right)\Delta_\mathrm{d}f\\
			& \quad -\Int\chi'\left(\mid\!f^i f^k\nabla_i f_k\!\mid^2-\mid\!f^i\bar f^{j}f_{i\bar j}\!\mid^2\right)\\
			& \quad -\Int\chi\,\operatorname{Im}\left(f^{i}T_{i}\right)\operatorname{Im}\left(f^if^k\nabla_i f_k\right).
		\end{split}
	\end{equation}
\end{lemma}
\begin{proof}
	By definition,
	\[
	2\operatorname{Re}\left(f^if^k\bar{f}^{j}\left(\nabla_i f_{k\bar j}-\nabla_{\bar j}\nabla_i f_k\right)\right)
	=-2\,R_{i\bar j k\bar l}f^i\bar{f}^{j}f^k\bar{f}^{l}.
	\]
	Integration by parts yields
	\[
	\begin{split}
		\Int\chi\,2\operatorname{Re}\left(f^if^k\bar f^{j}\nabla_i f_{k\bar j}\right)
		& = \Int\chi\left(f^i\bar f^{j}f_{i\bar j}\,\Delta_{\mathrm{d}}f
		-2\mid\!f^{k}f_{k\bar j}\!\mid^2\right)\\
		& \quad -\Int\chi\,2\operatorname{Re}\left(g^{i\bar j}f^k f^p f_{p\bar j}\nabla_k f_i\right)\\
		& \quad -\Int f^i\bar f^{j}f_{i\bar j}\,2\operatorname{Re}\left(f^k\partial_{k}\chi\right),
	\end{split}
	\]
	and
	\[
	\begin{split}
		-\Int\chi\,2\operatorname{Re}\left(f^if^k\bar f^{j}\nabla_{\bar j}\nabla_i f_k\right)
		& = \Int\chi\left(
		2\mid\!f^k\nabla_k f_i\!\mid^2
		+2\operatorname{Re}\left(g^{i\bar j}f^k\bar f^{l}
		\nabla_k f_i\nabla_{\bar j}f_{\bar l}\right)\right)\\
		& \quad +\Int\chi\,2\operatorname{Re}\left(\left(\Delta f+\bar f^{j}\bar T_j\right)f^if^k\nabla_i f_k\right)\\
		& \quad +\Int2\operatorname{Re}\left(f^if^k\nabla_i f_k\bar f^{l}\partial_{\bar l}\chi\right)
	\end{split}
	\]
	Since
	\[
	-f^i\bar f^{j}f_{i\bar j}\,2\operatorname{Re}\left(f^k\partial_{k}\chi\right)
	+2\operatorname{Re}\left(f^if^k\nabla_i f_k\bar f^{l}\partial_{\bar l}\chi\right)
	=2\chi'\left(\mid\!f^i f^k\nabla_i f_k\!\mid^2-\mid\!f^i\bar f^{j}f_{i\bar j}\!\mid^2\right)
	\]
	and
	\[
	-\operatorname{Re}\left(\bar{f}^{j}\bar{T}_{j}f^if^k\nabla_i f_k\right)
	+\frac{1}{2}\operatorname{Re}\left(f^if^k\nabla_i f_k\right)\Phi(f)
	=-\operatorname{Im}\left(f^{i}T_{i}\right)\operatorname{Im}\left(f^if^k\nabla_i f_k\right),
	\]
	combining all the above identities yields \eqref{eq:generalweight}.
\end{proof}

\begin{lemma} \label{thm:hermitian}
	Let $(X, \omega)$ be a compact Hermitian surface. For every $f\in\mathcal{A}^{0}_{\mathbb{R}}(X)$ and every $t\in\R$, the following identity holds
	\begin{equation}\label{eq:intro-hermitian}
		\begin{split}
			\Int F^{t}_{i\bar j k\bar l}f^i\bar f^{j}f^k\bar f^{l}
			& = -\Int\left\lvert f^k\nabla_if_k-\frac12\left(\Delta_\mathrm{d}f\right)f_i\right\rvert^2
			-\Int\left(\mid\!f^k\nabla_k f_i\!\mid^2+\mid\!\nabla f\!\mid^2\mid\!f_{i\bar j}\!\mid^2\right)\\
			& \quad +\frac{1}{2}\Int\mid\!\nabla f\!\mid^2\left(\Delta_{\mathrm{d}}f\right)^2+2\Int f^i\bar f^{j}f_{i\bar j}\,\Delta f\\
			& \quad +\Int\mid\!\nabla f\!\mid^2\,\operatorname{Re}\left(g^{i\bar j}f^{k}
			\left(T_i f_{k\bar j}-\bar{T}_{j}\nabla_{i}f_{k}\right)\right)\\
			& \quad -\Int\operatorname{Re}\left(f^if^k\nabla_i f_k\right)\,\Delta_\mathrm{d}f\\
			& \quad 
			+\frac12\Int\mid\!\nabla f\!\mid^2\Phi(f)\, \Delta_\mathrm{d}f
			+\frac14\Int\mid\!\nabla f\!\mid^2\left(\Phi(f)\right)^2\\
			& \quad +t^2\Int\left(\mid\!\nabla f\!\mid^2\mid\!f^iT_i\!\mid^2
			-\mid\!\tau\!\mid^2\mid\!\nabla f\!\mid^4\right).
		\end{split}
	\end{equation}
\end{lemma}
\begin{proof}
	By \eqref{eq:bar-hessian-commute},
	\[
	\begin{split}
		2\operatorname{Re}\left(g^{i\bar{j}}f^{k}\bar{f}^{l}\nabla_{k}f_{i}\nabla_{\bar{j}}f_{\bar{l}}\right)
		& = 2\mid\!f^k\nabla_i f_k\!\mid^2+\mid\!\nabla f\!\mid^2 2\operatorname{Re}\left(g^{i\bar j}\bar f^{l}T_i\nabla_{\bar j}f_{\bar l}\right)\\
		& \quad -2\operatorname{Re}\left(\bar{f}^{j}\bar{T}_{j}f^{k}f^{i}\nabla_{k}f_{i}\right).
	\end{split}
	\]
	By the above identity and \eqref{eq:bar-hessian-commute}, taking $\chi=1$ in \eqref{eq:generalweight} gives
	\begin{equation}\label{eq:quartic-commute}
		\begin{split}
			\Int R_{i\bar j k\bar l}f^i\bar f^{j}f^k\bar f^{l}
			& = -\Int\left(\mid\!f^k\nabla_k f_i\!\mid^2+\mid\!f^k\nabla_i f_k\!\mid^2\right)
			+\Int\mid\!f^{k}f_{k\bar j}\!\mid^2\\
			& \quad -\frac{1}{2}\Int f^i\bar f^{j}f_{i\bar j}\,\Delta_\mathrm{d}f
			-\Int\operatorname{Re}\left(f^if^k\nabla_i f_k\right)\Delta f\\
			& \quad +\Int\operatorname{Re}\left(g^{i\bar j}f^k f^p f_{p\bar j}\nabla_k f_i\right)
			-\Int\mid\!\nabla f\!\mid^2\,\operatorname{Re}\left(g^{i\bar j}\bar f^{l}T_i\nabla_{\bar j}f_{\bar l}\right).
		\end{split}
	\end{equation}
	
	By \eqref{eq:bar-hessian-commute},
	\[
	\begin{split}
		2\Real\left(g^{i\bar j}f^k f^p f_{p\bar j}\nabla_k f_i\right)
		& = 2\Real\left(g^{i\bar j}f^k f^p f_{p\bar j}\nabla_i f_k\right)
		-f^i\bar f^{j}f_{i\bar j}\,\Phi(f)\\
		& \quad +\mid\!\nabla f\!\mid^2\,2\Real\left(g^{i\bar j}f^kT_if_{k\bar j}\right).
	\end{split}
	\]
	Integration by parts gives (using \eqref{eq:bar-hessian-commute})
	\[
	\begin{split}
		\Int2\Real\left(g^{i\bar j}f^k f^p f_{p\bar j}\nabla_i f_k\right)
		& = -\Int\mid\!\nabla f\!\mid^2\,2\Real\left(g^{i\bar j}f^{p}\left(\nabla_i f_{p\bar j}+T_i f_{p\bar j}\right)\right)\\
		& \quad -2\Int\mid\!\nabla f\!\mid^2\mid\!f_{i\bar j}\!\mid^2
		-2\Int\mid\!f^{k}f_{k\bar j}\!\mid^2\\
		& = -\Int\mid\!\nabla f\!\mid^2\, 2\Real\left(f^k\partial_k\Delta f\right)
		-\Int\mid\!\nabla f\!\mid^2\Phi(f)\, \Delta f\\
		& \quad -2\Int\mid\!\nabla f\!\mid^2\mid\!f_{i\bar j}\!\mid^2
		-2\Int\mid\!f^{k}f_{k\bar j}\!\mid^2,
	\end{split}
	\]
	which implies
	\begin{equation} \label{eq:mixed}
		\begin{split}
			\Int2\Real\left(g^{i\bar j}f^k f^p f_{p\bar j}\nabla_k f_i\right)
			& = -\Int\mid\!\nabla f\!\mid^2\, 2\Real\left(f^k\partial_k\Delta f\right)
			-\Int\mid\!\nabla f\!\mid^2\Phi(f)\, \Delta f\\
			& \quad -2\Int\mid\!\nabla f\!\mid^2\mid\!f_{i\bar j}\!\mid^2
			-2\Int\mid\!f^{k}f_{k\bar j}\!\mid^2\\
			& \quad -\Int f^i\bar f^{j}f_{i\bar j}\,\Phi(f)
			+\Int\mid\!\nabla f\!\mid^2\,2\Real\left(g^{i\bar j}f^kT_if_{k\bar j}\right).
		\end{split}
	\end{equation}
	Integration by parts yields (using \eqref{eq:bar-hessian-commute})
	\[
	\begin{split}
		-\Int\mid\!\nabla f\!\mid^2\, 2\Real\left(f^k\partial_k\Delta f\right)
		& = -\Int\mid\!\nabla f\!\mid^2\left(\Delta_{\mathrm{d}}f\right)\left(\Delta f\right)
		+\Int2\Real\left(f^k\partial_k\mid\!\nabla f\!\mid^2\right)\,\Delta f\\
		& = 2\Int\mid\!\nabla f\!\mid^2\left(\Delta f\right)^2+\Int\mid\!\nabla f\!\mid^2\Phi(f)\, \Delta f\\
		& \quad +2\Int f^i\bar f^{j}f_{i\bar j}\,\Delta f
		+\Int2\Real\left(f^if^k\nabla_i f_k\right)\,\Delta f.
	\end{split}
	\]
	Substitute the above equation into \eqref{eq:mixed}, and then substitute the resulting equation into \eqref{eq:quartic-commute} to obtain
	\[
	\begin{split}
		\Int R_{i\bar j k\bar l}f^i\bar f^{j}f^k\bar f^{l}
		& = -\Int\left(\mid\!f^k\nabla_k f_i\!\mid^2+\mid\!f^k\nabla_i f_k\!\mid^2
		+\mid\!\nabla f\!\mid^2\mid\!f_{i\bar j}\!\mid^2\right)\\
		& \quad +\Int\mid\!\nabla f\!\mid^2\left(\Delta f\right)^2+2\Int f^i\bar f^{j}f_{i\bar j}\,\Delta f\\
		& \quad +\Int\mid\!\nabla f\!\mid^2\,\Real\left(g^{i\bar j}f^{k}
		\left(T_i f_{k\bar j}-\bar{T}_{j}\nabla_{i}f_{k}\right)\right).
	\end{split}
	\]
	Since
	\[
	\left\lvert f^k\nabla_i f_k-\frac12\left(\Delta_\mathrm{d}f\right)f_i\right\rvert^2
	=\mid\!f^k\nabla_i f_k\!\mid^2
	-\Real\left(f^if^k\nabla_i f_k\right)\,\Delta_\mathrm{d}f
	+\frac14\mid\!\nabla f\!\mid^2\left(\Delta_\mathrm{d}f\right)^2,
	\]
	combining the above two identities yields
	\[
	\begin{split}
		\Int R_{i\bar j k\bar l}f^i\bar f^{j}f^k\bar f^{l}
		& = -\Int\left\lvert f^k\nabla_if_k-\frac12\left(\Delta_\mathrm{d}f\right)f_i\right\rvert^2
		-\Int\left(\mid\!f^k\nabla_k f_i\!\mid^2+\mid\!\nabla f\!\mid^2\mid\!f_{i\bar j}\!\mid^2\right)\\
		& \quad +\frac{1}{2}\Int\mid\!\nabla f\!\mid^2\left(\Delta_{\mathrm{d}}f\right)^2+2\Int f^i\bar f^{j}f_{i\bar j}\,\Delta f\\
		& \quad +\Int\mid\!\nabla f\!\mid^2\,\Real\left(g^{i\bar j}f^{k}
		\left(T_i f_{k\bar j}-\bar{T}_{j}\nabla_{i}f_{k}\right)\right)\\
		& \quad -\Int\Real\left(f^if^k\nabla_i f_k\right)\,\Delta_\mathrm{d}f\\
		& \quad 
		+\frac12\Int\mid\!\nabla f\!\mid^2\Phi(f)\, \Delta_\mathrm{d}f
		+\frac14\Int\mid\!\nabla f\!\mid^2\left(\Phi(f)\right)^2.
	\end{split}
	\]
	Then \eqref{eq:intro-hermitian} follows from \eqref{eq:rankone}.
\end{proof}

Since
\[
\begin{split}
	\Int\mid\!\nabla f\!\mid^2\Real\left(f^k\partial_k\Delta_\mathrm{d}f\right)
	& = \frac{1}{2}\Int\mid\!\nabla f\!\mid^2\left(\Delta_\mathrm{d}f\right)^2
	-\Int\Real\left(f^if^k\nabla_i f_k\right)\,\Delta_\mathrm{d}f\\
	& \quad -\Int f^i\bar f^{j}f_{i\bar j}\,\Delta_\mathrm{d}f,
\end{split}
\]
we obtain from \eqref{eq:intro-hermitian} that
\begin{equation}\label{11031}
	\begin{split}
		\Int\mid\!\nabla f\!\mid^2\Real\left(f^k\partial_k\Delta_\mathrm{d}f\right)
		& = \Int F^{t}_{i\bar j k\bar l}f^i\bar f^{j}f^k\bar f^{l}+\Int\left\lvert f^k\nabla_if_k-\frac12\left(\Delta_\mathrm{d}f\right)f_i\right\rvert^2\\
		& \quad +\Int\left(\mid\!f^k\nabla_k f_i\!\mid^2+\mid\!\nabla f\!\mid^2\mid\!f_{i\bar j}\!\mid^2\right)+\Int f^i\bar f^{j}f_{i\bar j}\,\Phi(f)\\
		& \quad -\Int\mid\!\nabla f\!\mid^2\,\operatorname{Re}\left(g^{i\bar j}f^{k}
		\left(T_i f_{k\bar j}-\bar{T}_{j}\nabla_{i}f_{k}\right)\right)\\
		& \quad 
		-\frac12\Int\mid\!\nabla f\!\mid^2\Phi(f)\, \Delta_\mathrm{d}f
		-\frac14\Int\mid\!\nabla f\!\mid^2\left(\Phi(f)\right)^2\\
		& \quad -t^2\Int\left(\mid\!\nabla f\!\mid^2\mid\!f^iT_i\!\mid^2
		-\mid\!\tau\!\mid^2\mid\!\nabla f\!\mid^4\right).
	\end{split}
\end{equation}
When $\tau=0$, we get
\[
\begin{split}
	\Int\mid\!\nabla f\!\mid^2\operatorname{Re}\left(f^k\partial_k\Delta_\mathrm{d}f\right)
	& = \Int R_{i\bar j k\bar l}f^i\bar f^{j}f^k\bar f^{l}
	+\Int\left\lvert f^k\nabla_if_k-\frac12\left(\Delta_\mathrm{d}f\right)f_i\right\rvert^2\\
	& \quad +\Int\left(\mid\!f^k\nabla_k f_i\!\mid^2+\mid\!\nabla f\!\mid^2\mid\!f_{i\bar j}\!\mid^2\right).
\end{split}
\]
\begin{remark}\label{11071}
	If $\Delta_\mathrm{d}f=\psi(f)$, then
	\[
	\Int\mid\!\nabla f\!\mid^2\operatorname{Re}\left(f^k\partial_k\Delta_\mathrm{d}f\right)
	=\Int\psi'(f)\mid\!\nabla f\!\mid^4.
	\]
	This applies both to eigenfunctions and to solutions of nonlinear Euler-Lagrange equations.
	
	If $\omega$ is K\"ahler and $\Delta_{\bar\partial}f=\mu f$ for some $\mu\in\mathbb R$, then $\tau=0$ and $\Delta_{\mathrm d}f=2\mu f$. With these substitutions, \eqref{11031} recovers the quartic eigenfunction identity of Wang and Yang in complex dimension two.
\end{remark}

To handlethe torsion terms, we introduce the following term. For every $f\in\mathcal{A}^{0}_{\mathbb{R}}(X)$, define
\begin{equation}\label{eq:Sf-invariant}
	\begin{split}
		S(f)
		&=
		\frac12
		\left\lvert f^k\nabla_i f_k+\bar f^l f_{i\bar l}
		-\mid\!\nabla f\!\mid^2T_i+
		\left(\frac12\Phi(f)+\frac{2}{3}\sqrt{-1}\operatorname{Im}\left(f^kT_k\right)
		\right)f_i\right\rvert^2\\
		& \quad +
		\left\lvert f^k\nabla_i f_k-\bar f^l f_{i\bar l}
		+\mid\!\nabla f\!\mid^2T_i
		-\left(\frac12\Delta_\mathrm{d}f+\frac34\Phi(f)
		+\frac{4}{3}\sqrt{-1}\operatorname{Im}\left(f^kT_k\right)
		\right)f_i\right\rvert^2\\
		& \quad +\frac14 \mid\!\nabla f\!\mid^2
		\left(\Delta_\mathrm{d}f-\frac12 \Phi(f)\right)^2.
	\end{split}
\end{equation}
Here each expression inside the first two norms is a $(1,0)$-form. All derivatives are Chern covariant derivatives. Thus \eqref{eq:Sf-invariant} defines a smooth nonnegative function on all of $X$, independently of any choice of local frames. Moreover, on $\left\{\mathrm{d}f=0\right\}\subseteq X$, we have $S(f)=0$.

\begin{lemma} \label{thm:quartic}
	Let $(X, \omega)$ be a compact Gauduchon surface. For every $f\in\mathcal{A}^{0}_{\mathbb{R}}(X)$ and every $t\in\R$, the following identity holds
	\begin{equation} \label{eq:quartic-remainder}
		\begin{split}
			\Int\mid\!\nabla f\!\mid^2\operatorname{Re}\left(f^k\partial_k\Delta_\mathrm{d}f\right)
			& = \Int F^t_{i\bar j k\bar l}f^i\bar f^{j}f^k\bar f^{l}
			+\frac12\Int\left(\left\lvert\nabla\mid\!\nabla f\!\mid^2\right\rvert^2-\mid\!\tau\!\mid^2\mid\!\nabla f\!\mid^4\right)\\
			& \quad +t^2\Int\left(\mid\!\tau\!\mid^2\mid\!\nabla f\!\mid^4-\mid\!\nabla f\!\mid^2\mid\!f^iT_i\!\mid^2\right)\\
			& \quad +\frac13\Int\mid\!\nabla f\!\mid^2\left(\Imag\left(f^iT_i\right)\right)^2+\Int S(f).
		\end{split}
	\end{equation}
\end{lemma}
\begin{proof}
	By \eqref{eq:bar-hessian-commute},
	\[
	\partial_i\mid\!\nabla f\!\mid^2
	=f^k\nabla_i f_k+\bar f^j f_{i\bar j},\quad
	f^k\nabla_k f_i
	=f^k\nabla_i f_k+\mid\!\nabla f\!\mid^2T_i-\left(f^kT_k\right)f_i.
	\]
	Since
	\[
	\begin{split}
		\left\lvert\partial f\wedge\sqrt{-1}\,\partial\bar\partial f\right\rvert^{2}\,\frac{\omega^2}{2!}
		& = \left(\Delta f\right)\,\sqrt{-1}\,\partial f \wedge \bar\partial f \wedge\sqrt{-1}\,\partial\bar\partial f\\
		& \quad -\frac{1}{2}\mid\!\nabla f\!\mid^{2}\,\sqrt{-1}\,\partial\bar\partial f\wedge\sqrt{-1}\,\partial\bar\partial f,
	\end{split}
	\]
	we have
	\begin{equation}\label{11131}
		\mid\!\nabla f\!\mid^2\mid\!f_{i\bar j}\!\mid^2-2\mid\!f^k f_{k\bar j}\!\mid^2
		=
		\mid\!\nabla f\!\mid^2\left(\Delta f\right)^2-2f^i\bar f^j f_{i\bar j}\Delta f.
	\end{equation}
	Expanding the three squares in \eqref{eq:Sf-invariant} and using these identities gives
	\begin{align}
		&\mathrel{\phantom{=}}
		\frac12\left(\left\lvert\nabla\mid\!\nabla f\!\mid^2\right\rvert^2-\mid\!\tau\!\mid^2\mid\!\nabla f\!\mid^4
		\right)+\frac13\mid\!\nabla f\!\mid^2\left(\Imag\left(f^kT_k\right)\right)^2+S(f) \notag\\
		& = \left\lvert f^k\nabla_i f_k-\frac12\left(\Delta_{\mathrm d}f\right)f_i\right\rvert^2
		+\mid\!f^k\nabla_k f_i\!\mid^2+\mid\!\nabla f\!\mid^2\mid\!f_{i\bar j}\!\mid^2
		+f^i\bar f^j f_{i\bar j}\,\Phi(f) \label{11032}\\
		& \quad -\frac12 \mid\!\nabla f\!\mid^2\Phi(f)\Delta_{\mathrm d}f
		-\frac14 \mid\!\nabla f\!\mid^2\left(\Phi(f)\right)^2
		-\mid\!\nabla f\!\mid^2\Phi\left(\mid\!\nabla f\!\mid^2\right) \notag\\
		& \quad -\mid\!\nabla f\!\mid^2\operatorname{Re}\left(g^{i\bar j}f^k\left(T_i f_{k\bar j}-\bar T_j\nabla_i f_k\right)\right).\notag
	\end{align}
	By Lemma \ref{11041},
	\[
	\Int\mid\!\nabla f\!\mid^2\Phi(\mid\!\nabla f\!\mid^2)
	=\frac12\Int\Phi(\mid\!\nabla f\!\mid^4)=0.
	\]
	Integrating \eqref{11032} and substituting it into \eqref{11031}, we obtain \eqref{eq:quartic-remainder}.
\end{proof}

The derivative term retained in \eqref{eq:quartic-remainder} gives a criterion that does not require pointwise positive sectional curvature.

For every non-zero $u\in W^{1,2}(X;\mathbb R)$, define
\[
\mathcal{E}_\omega(u)
=\Int\left(\frac12\mid\!\nabla u\!\mid^2+\mathcal{H}_{\omega}(x)\, u^2\right),
\]
where $\mathcal H_\omega(x)$ is defined by \eqref{eq:effective-quartic-potential}. Then
\[
\sigma_\omega=\inf_{u}\frac{\mathcal{E}_\omega(u)}{\Int u^2},
\]
where $\sigma_\omega$ is defined by \eqref{eq:curvature-form-bottom}.
\begin{theorem}\label{11062}
	Let $(X,\omega)$ be a compact connected Gauduchon surface, and let $\sigma_\omega$ be defined by \eqref{eq:curvature-form-bottom}. Then
	\[
	\lambda_1>\sigma_\omega.
	\]
	More generally, let $f\in\mathcal A^0_{\mathbb R}(X)$ be a nonconstant solution of $\Delta_{\mathrm d}f=\psi(f)$, where $\psi$ is $C^{1}$ on a neighborhood of $f(X)$. Then
	\[
	\frac{\Int\psi'(f)\mid\!\nabla f\!\mid^4}{\Int\mid\!\nabla f\!\mid^4}>\sigma_\omega.
	\]
	Consequently, if $\psi'(y)\leq\sigma_\omega$ for every $y\in f(X)$, then $f$ is a constant.
\end{theorem}

\begin{proof}
	Let $f$ be a nonconstant smooth real solution of $\Delta_{\mathrm d}f=\psi(f)$. By \eqref{eq:quartic-remainder} with $t=0$ and the definition of $\mathcal{H}_{\omega}$,
	\[
	\begin{split}
		\Int\psi'(f)\mid\!\nabla f\!\mid^4
		&\geq
		\Int\left(
		\frac12\left\lvert\nabla\mid\!\nabla f\!\mid^2\right\rvert^2
		+\mathcal{H}_{\omega}\mid\!\nabla f\!\mid^4\right)\\
		&=\mathcal{E}_\omega(\mid\!\nabla f\!\mid^2)
		\geq\sigma_\omega\Int\mid\!\nabla f\!\mid^4.
	\end{split}
	\]
	The last inequality follows from \eqref{eq:curvature-form-bottom} by choosing the test function $u=\mid\!\nabla f\!\mid^2$.
	
	The first inequality is strict. Indeed, equality would imply $S(f)\equiv0$. The last square in \eqref{eq:Sf-invariant} would then give
	\[
	\mid\!\nabla f\!\mid^2\left(\Delta_\mathrm{d}f-\frac12 \Phi(f)\right)^2=0.
	\]
	Thus $\Delta_{\mathrm d}f=\frac12\Phi(f)$ wherever $\mathrm df\ne0$. The same identity holds on the interior of the critical set, where $f$ is locally constant, and on its boundary by continuity. Multiplying by $f$ and integrating, we obtain from Lemma \ref{11041} that
	\[
	\Int\mid\!\mathrm df\!\mid^2=\Int f\,\Delta_{\mathrm d}f=\frac12\Int f\,\Phi(f)=\frac14\Int\Phi(f^2)=0,
	\]
	This contradicts the assumption that $f$ is nonconstant. Since $\Int\mid\!\nabla f\!\mid^4>0$, taking a first eigenfunction and $\psi(y)=\lambda_1\,y$ gives the eigenvalue estimate. Finally, if $\psi'\leq\sigma_\omega$ on $f(X)$, the strict inequality rules out every nonconstant solution.
\end{proof}

\begin{lemma}
	Let $(X,\omega)$ be a compact Hermitian surface. For every $f\in\mathcal A^0_{\mathbb R}(X)$ and every $t\in\mathbb R$, the following identity holds
	\begin{equation}\label{11091}
		\begin{split}
			0
			& = \int_{X_+}\frac{1}{\mid\!\nabla f\!\mid^2}F^t_{i\bar j k\bar l}
			f^i\bar f^j f^k\bar f^l
			+\int_{X_+}\frac{1}{\mid\!\nabla f\!\mid^2}\left(\mid\!f^k\nabla_k f_i\!\mid^2
			-\mid\!f^kf_{k\bar{j}}\!\mid^2\right)\\
			& \quad +\int_{X_+}\frac{1}{\mid\!\nabla f\!\mid^2}\operatorname{Re}\left(g^{i\bar j}f^k\nabla_k f_i\left(\bar f^l\nabla_{\bar j}f_{\bar l}-f^p f_{p\bar j}\right)\right)\\
			& \quad +\int_{X_+}\frac{1}{\mid\!\nabla f\!\mid^4}\left(\left(f^i\bar f^j f_{i\bar j}\right)^2-\mid\!f^if^k\nabla_i f_k\!\mid^2\right)\\
			& \quad +\frac12\int_{X_+}\frac{1}{\mid\!\nabla f\!\mid^2}\left(f^i\bar f^j f_{i\bar j}-\operatorname{Re}\left(f^if^k\nabla_i f_k\right)\right)\, \Delta_{\mathrm d}f\\
			& \quad +\int_{X_+}\frac{1}{\mid\!\nabla f\!\mid^2}\Imag\left(f^iT_i\right)\,\Imag\left(f^if^k\nabla_i f_k\right)
			+t^2\int_{X_+}\left(\mid\!\tau\!\mid^2\mid\!\nabla f\!\mid^2-\mid\!f^iT_i\!\mid^2\right).
		\end{split}
	\end{equation}
\end{lemma}
\begin{proof}
	For $\varepsilon>0$, apply \eqref{eq:generalweight} with
	\[
	\chi_\varepsilon(\mid\!\nabla f\!\mid^2)=\frac{1}{\mid\!\nabla f\!\mid^2+\varepsilon}.
	\]
	Then $\chi_\varepsilon$ and $\chi'_\varepsilon$ are smooth weights on all of $X$.
	
	Consider the case: $\varepsilon\rightarrow0$. By the Cauchy-Schwarz inequality,
	\[
	\frac{1}{\left(\mid\!\nabla f\!\mid^2+\varepsilon\right)^2}\left(\mid\!f^if^k\nabla_i f_k\!\mid^2+\mid\!f^i\bar f^j f_{i\bar j}\!\mid^2\right)
	\leq\mid\!\nabla_i f_k\!\mid^2+\mid\!f_{i\bar j}\!\mid^2.
	\]
	Together with Young's inequality, these estimates show that all regularized integrands are bounded in absolute value by
	\[
	C\left(\mid\!R\!\mid\mid\!\nabla f\!\mid^2+\mid\!\nabla_i f_k\!\mid^2
	+\mid\!f_{i\bar j}\!\mid^2+\left(\Delta_{\mathrm d}f\right)^2+\mid\!\tau\!\mid^2\mid\!\nabla f\!\mid^2\right),
	\]
	where $C$ is independent of $\varepsilon$. This function is integrable because $X$ is compact.
	
	On $\{\mid\!\nabla f\!\mid^2=0\}$, every regularized integrand vanishes. Dominated convergence therefore gives an identity over $X_+$, with each quotient extended by zero on the critical set. Moving all terms to the left gives the asserted formula for $t=0$. In particular, the term involving $\chi_\varepsilon'$ becomes
	\[
	\int_{X_+}
	\frac{1}{\mid\!\nabla f\!\mid^4}
	\left(\left(f^i\bar f^j f_{i\bar j}\right)^2-\mid\!f^if^k\nabla_i f_k\!\mid^2\right).
	\]
	While on $X_+$, \eqref{eq:rankone} gives 
	\[
	\frac{1}{\mid\!\nabla f\!\mid^2}R_{i\bar j k\bar l}f^i\bar f^j f^k\bar f^l
	=
	\frac{1}{\mid\!\nabla f\!\mid^2}F^t_{i\bar j k\bar l}f^i\bar f^j f^k\bar f^l
	+t^2\left(\mid\!\tau\!\mid^2\mid\!\nabla f\!\mid^2-\mid\!f^iT_i\!\mid^2\right).
	\]
	Substitution proves the formula for every real $t$.
\end{proof}

Define the nonnegative remainder on $X_+$ by
\begin{equation}\label{eq:Qstar}
	\begin{split}
		Q(f)
		& =
		\left\lvert\frac{1}{\mid\!\nabla f\!\mid^2}f^if^k\nabla_i f_k
		-\frac{1}{2\mid\!\nabla f\!\mid^2}f^i\bar f^j f_{i\bar j}-\frac14\Delta_{\mathrm d}f
		+\frac1{20}\Phi(f)+\frac35\sqrt{-1}\,\Imag\left(f^iT_i\right)\right\rvert^2\\
		& \quad +\frac{2}{\mid\!\nabla f\!\mid^2}\left\lvert f^k\nabla_i f_k-\frac14\bar f^j f_{i\bar j}+\frac45\left(\mid\!\nabla f\!\mid^2T_i-\left(f^kT_k\right)f_i\right)
		-\frac{1}{\mid\!\nabla f\!\mid^2}\left(f^kf^p\nabla_p f_k-\frac14f^k\bar f^l f_{k\bar l}\right)f_i\right\rvert^2\\
		&+\frac78
		\left(\frac{1}{\mid\!\nabla f\!\mid^2}f^i\bar f^j f_{i\bar j}+\frac{13}{28}\Delta_{\mathrm d}f+\frac{57}{140}\Phi(f)\right)^2.
	\end{split}
\end{equation}
On $X_+$, the Cauchy-Schwarz inequality gives
\[
0\leq Q(f)\leq C\left(\mid\!\nabla_i f_k\!\mid^2+\mid\!f_{i\bar j}\!\mid^2+\left(\Delta_{\mathrm d}f\right)^2+\mid\!\tau\!\mid^2\mid\!\nabla f\!\mid^2\right).
\]
We extend $Q(f)$ by zero to $X\setminus X_+$. The resulting function is integrable.
\begin{lemma} \label{thm:completion}
	Let $(X,\omega)$ be a compact Gauduchon surface, and let $f\in\mathcal A^0_{\mathbb R}(X)$ satisfy $\Delta_{\mathrm d}f=\psi(f)$, where $\psi$ is $C^{1}$ on a neighborhood of $f(X)$. For every $t\in\mathbb R$, the following identity holds
	\begin{align}
		&\mathrel{\phantom{=}}
		\frac{225}{448}\Int\psi'(f)\mid\!\nabla f\!\mid^2
		+\frac1{25}\Int\left(7\mid\!\tau\!\mid^2\mid\!\nabla f\!\mid^2+8\mid\!f^iT_i\!\mid^2\right)
		-\int_{X_+}\frac{1}{\mid\!\nabla f\!\mid^2}F^t_{i\bar j k\bar l}f^i\bar f^j f^k\bar f^l\label{eq:new-wy-exact}\\
		& =\Int Q(f)+\frac{11}{1120}\Int\left(\Real\left(f^iT_i\right)\right)^2
		+\frac6{25}\Int\left(\Imag\left(f^iT_i\right)\right)^2
		+t^2\Int\left(\mid\!\tau\!\mid^2\mid\!\nabla f\!\mid^2-\mid\!f^iT_i\!\mid^2\right).\notag
	\end{align}
\end{lemma}
\begin{proof}
	Let $J_{f,t}$ denote the sum of the noncurvature integrands in \eqref{11091}, so that
	\[
	\int_{X_+}J_{f,t}
	=
	-\int_{X_+}
	\frac{1}{\mid\!\nabla f\!\mid^2}F^t_{i\bar j k\bar l}f^i\bar f^j f^k\bar f^l.
	\]
	
	To expand the middle square in \eqref{eq:Qstar}, we use the orthogonal projection formula: on $X_+$, for a $(1,0)$-form $\eta$, 
	\[
	\left\lvert\eta_i-\frac{1}{\mid\!\nabla f\!\mid^2}\left(f^k\eta_k\right)f_i\right\rvert^2
	=
	\mid\!\eta\!\mid^2-\frac{1}{\mid\!\nabla f\!\mid^2}\mid\!f^k\eta_k\!\mid^2.
	\]
	
	Expanding the squares in \eqref{eq:Qstar} and simplifying by means of \eqref{eq:bar-hessian-commute}, \eqref{11131}, and the preceding identities, we obtain on $X_+$:
	\begin{align*}
		&\mathrel{\phantom{=}}
		Q(f)
		+\frac{11}{1120}\left(\Real\left(f^iT_i\right)\right)^2
		+\frac6{25}\left(\Imag\left(f^iT_i\right)\right)^2
		+t^2\left(\mid\!\tau\!\mid^2\mid\!\nabla f\!\mid^2-\mid\!f^iT_i\!\mid^2\right)\\
		& = J_{f,t}
		+\frac{225}{896}\left(\Delta_{\mathrm d}f\right)^2
		+\frac1{25}\left(7\mid\!\tau\!\mid^2\mid\!\nabla f\!\mid^2+8\mid\!f^iT_i\!\mid^2\right)+\frac9{16}\left(\mid\!f_{i\bar j}\!\mid^2-\left(\Delta f\right)^2\right)\\
		& \quad +\frac1{10}\Phi\left(\mid\!\nabla f\!\mid^2\right)
		+\frac{137}{448}\left(\Delta_{\mathrm d}f\right)\Phi(f).
	\end{align*}
	
	We extend $Q(f)$ by zero to the critical set. Both Chern Hessians, and hence $\Delta_{\mathrm d}f$, vanish almost everywhere on $X\setminus X_+$, as explained after
	Lemma \ref{11041}. Also, $\mid\!\nabla f\!\mid^2$, $\nabla\mid\!\nabla f\!\mid^2$, and $f^iT_i$ vanish there pointwise. We may therefore integrate the displayed identity over $X_+$ and take the integrals of its correction terms over all of $X$.
	
	The last three terms have zero integral by Lemma \ref{11041}. In particular, if $A'=\psi$, then
	\[
	\Int\left(\Delta_{\mathrm d}f\right)\Phi(f)
	=\Int\Phi\left(A(f)\right)=0.
	\]
	Finally, integration by parts gives
	\[
	\Int\left(\Delta_{\mathrm d}f\right)^2=\Int\psi(f)\,\Delta_{\mathrm d}f
	=2\Int\psi'(f)\mid\!\nabla f\!\mid^2.
	\]
	Substitution proves \eqref{eq:new-wy-exact}. Since
	\[
	Q(f)\geq0, \quad \mid\!\tau\!\mid^2\mid\!\nabla f\!\mid^2\geq\mid\!f^iT_i\!\mid^2,
	\]
	the right-hand side of \eqref{eq:new-wy-exact} is nonnegative.
\end{proof}

The Hermitian inner product on $(2,0)$-forms gives the pointwise identity
\begin{equation}\label{11075}
	\mid\!\tau\wedge\partial f\!\mid^2
	=\mid\!\tau\!\mid^2\mid\!\nabla f\!\mid^2-\mid\!f^iT_i\!\mid^2.
\end{equation}
Thus, on $X_+$, the squared norms of the components of $\tau$ parallel and orthogonal to $\partial f$ are, respectively,
\[
\frac{1}{\mid\!\nabla f\!\mid^2}\mid\!f^iT_i\!\mid^2, \quad
\frac{1}{\mid\!\nabla f\!\mid^2}\mid\!\tau\wedge\partial f\!\mid^2.
\]
Under the hypotheses of Lemma \ref{thm:completion}, taking $t=0$ in \eqref{eq:new-wy-exact} therefore gives
\begin{equation}\label{eq:directional-defect}
	\begin{split}
		\frac{225}{448}\Int\psi'(f)\mid\!\nabla f\!\mid^2
		& = \int_{X_+}\frac{1}{\mid\!\nabla f\!\mid^2}R_{i\bar j k\bar l}f^i\bar f^j f^k\bar f^l+\Int Q(f)\\
		& \quad -\int_{X_+}\left(\frac35\mid\!f^iT_i\!\mid^2+\frac7{25}\mid\!\tau\wedge\partial f\!\mid^2\right)\\
		& \quad +\frac{11}{4480}\Int\left(\Phi(f)\right)^2
		+\frac6{25}\Int\left(\Imag\left(f^iT_i\right)\right)^2.
	\end{split}
\end{equation}

\begin{lemma}\label{lem:strict-remainder}
	Let $(X,\omega)$ be a compact connected Hermitian surface. Suppose that $f$ is nonconstant smooth real function satisfying $\Delta_{\mathrm d}f=\psi(f)$, where $\psi$ is $C^{1}$ on a neighborhood of $f(X)$. Then
	\[
	\int_{X_+}Q(f)>0.
	\]
\end{lemma}
\begin{proof}
	Suppose that the integral vanishes. Since $Q(f)$ is nonnegative and continuous on $X_+$, it vanishes throughout $X_+$. The last square in \eqref{eq:Qstar} and the real part of the expression in its first square give
	\[
	\begin{aligned}
		\frac{1}{\mid\!\nabla f\!\mid^2}f^i\bar f^j f_{i\bar j}
		& = -\frac{13}{28}\Delta_{\mathrm d}f-\frac{57}{140}\Phi(f),\\
		\frac{1}{\mid\!\nabla f\!\mid^2}\Real\left(f^if^k\nabla_i f_k\right)
		& = \frac{1}{2\mid\!\nabla f\!\mid^2}f^i\bar f^j f_{i\bar j}
		+\frac14\Delta_{\mathrm d}f-\frac1{20}\Phi(f),
	\end{aligned}
	\]
	respectively. Then we get
	\begin{equation}\label{eq:strict-gradient-hessian}
		\Real\left(f^k\partial_k\mid\!\nabla f\!\mid^2\right)
		=-\frac{1}{56}\mid\!\nabla f\!\mid^2\left(25\Delta_{\mathrm d}f+37\Phi(f)\right).
	\end{equation}
	Both sides vanish wherever $\mathrm df=0$, so \eqref{eq:strict-gradient-hessian} holds on all of $X$.
	
	Let $\nabla_{\mathbb R}f$ be the gradient with respect to the underlying Riemannian metric. Since
	\[
	\mid\!\nabla_{\mathbb R}f\!\mid^2=2\mid\!\nabla f\!\mid^2,
	\]
	we have
	\[
	\operatorname{Hess}_{\mathbb R}f\left(\nabla_{\mathbb R}f, \nabla_{\mathbb R}f\right)
	=\frac12\nabla_{\mathbb R}f\left(\mid\!\nabla_{\mathbb R}f\!\mid^2\right)
	=\nabla_{\mathbb R}f\left(\mid\!\nabla f\!\mid^2\right)
	=2\Real\left(f^k\partial_k\mid\!\nabla f\!\mid^2\right).
	\]
	Thus \eqref{eq:strict-gradient-hessian} is equivalent to the global identity
	\begin{equation}\label{11061}
		\operatorname{Hess}_{\mathbb R}f\left(\nabla_{\mathbb R}f, \nabla_{\mathbb R}f\right)
		=-\frac{1}{56}\left(25\Delta_{\mathrm d}f+37\Phi(f)\right)\mid\!\nabla_{\mathbb R}f\!\mid^2.
	\end{equation}
	
	Choose $p\in X$ with
	\[
	f(p)=M:=\max_X f.
	\]
	We claim that $\Delta_{\mathrm d}f(p)>0$. Indeed,
	\[
	\Delta_{\mathrm d}f(p)=\psi(M)\geq0
	\]
	by the sign convention for $\Delta_{\mathrm d}$. If $\psi(M)=0$, then
	\[
	w=M-f\geq0
	\]
	satisfies $\Delta_{\mathrm d}w=a(x)w$, where
	\[
	a(x)=\int_0^1\psi'\left(M-\theta\,w(x)\right)\,\mathrm{d}\theta.
	\]
	Then we obtain
	\[
	-\left(\Delta_{\mathrm d}+\max\{-a,0\}\right)\,w=-\max\{a,0\}\,w\leq0.
	\]
	The strong minimum principle at $p$ gives $w\equiv0$, contrary to the nonconstancy of $f$.
	
	Let $A$ be the real self-adjoint endomorphism associated with $\operatorname{Hess}_{\mathbb R}f(p)$. Fix $v\in T_pX$ and differentiate \eqref{11061} twice along the geodesic with initial velocity $v$. At $p$, we obtain
	\[
	\nabla_{\mathbb R}f(p)=0, \quad \nabla^{\rm LC}_v\left(\nabla_{\mathbb R}f\right)=Av.
	\]
	Consequently, the only surviving terms give
	\[
	\operatorname{Hess}_{\mathbb R}f(p)\left(Av,Av\right)
	=-\frac{1}{56}\left(25\Delta_{\mathrm d}f+37\Phi(f)\right)(p)\mid\!Av\!\mid^2.
	\]
	Since $\Phi(f)(p)=0$, every nonzero eigenvalue of $A$ therefore equals
	\[
	-\frac{1}{56}\left(25\Delta_{\mathrm d}f+37\Phi(f)\right)(p)
	=-\frac{25}{56}\Delta_{\mathrm d}f(p).
	\]
	If $m=\operatorname{rank}A$, taking the trace gives
	\[
	-\Delta_{\mathrm d}f(p)
	=\operatorname{tr}A
	=-m\,\frac{25}{56}\,\Delta_{\mathrm d}f(p).
	\]
	As $\Delta_{\mathrm d}f(p)>0$, this would imply $m=\frac{56}{25}$, contradicting the fact that $m$ is an integer.
\end{proof}

\begin{proof}[Proof of Theorem \ref{thm:directional-spectral-gap}]
	Let $f$ be a real first eigenfunction of $\Delta_{\mathrm d}$. Applying \eqref{eq:directional-curvature-bound} to the complex gradient $f^i\partial_i$, we obtain on $X_+$
	\[
	\frac{1}{\mid\!\nabla f\!\mid^2}R_{i\bar j k\bar l}f^i\bar f^j f^k\bar f^l
	-\frac35\mid\!f^iT_i\!\mid^2
	-\frac7{25}\mid\!\tau\wedge\partial f\!\mid^2\geq\kappa\mid\!\nabla f\!\mid^2,
	\]
	where we have used \eqref{11075}. Taking $\psi(y)=\lambda_1\, y$ in \eqref{eq:directional-defect} gives
	\[
	\left(\frac{225}{448}\,\lambda_1-\kappa\right)\Int\mid\!\nabla f\!\mid^2\geq\Int Q(f)>0.
	\]
	The strict inequality follows from Lemma \ref{lem:strict-remainder}, since $f$ is nonconstant. Then
	\[
	\lambda_1>\frac{448}{225}\,\kappa.
	\]
	Now let $f$ be any nonconstant smooth real solution of $\Delta_{\mathrm d}f=\psi(f)$. Integration by parts gives
	\[
	\frac{\Int\psi'(f)\mid\!\nabla f\!\mid^2}{\Int\mid\!\nabla f\!\mid^2}
	=
	\frac{\Int\left(\Delta_{\mathrm d}f\right)^2}{\Int\mid\!\mathrm{d}f\!\mid^2}
	\geq\lambda_1>\frac{448}{225}\,\kappa.
	\]
	The first inequality follows by expanding $f$ in Laplace eigenfunctions; its constant component contributes to neither integral.
\end{proof}

\begin{proof}[Proof of Corollary \ref{cor:isotropic}]
	By \eqref{eq:rankone}, for every $v\in\left(T^{1,0}X\backslash\{0\}\right)$,
	\[
	\operatorname{HSC}^{t}(v)
	=\operatorname{HSC}^{0}(v)-t^2\left(\mid\!\tau\!\mid^2-\frac{\mid\!\tau(v)\!\mid^2}{\mid\!v\!\mid^2}\right).
	\]
	Thus the Chern condition is the weakest of the uniform bounds in \eqref{eq:main-pinching}. In \eqref{eq:new-wy-exact} the change in sectional curvature is exactly compensated by the explicit $t^2$-term. This cancellation concerns the sectional contraction of $F^t$; the two Ricci traces in Definition \ref{10112} have different formulas. Under \eqref{eq:main-pinching} it follows that
	\[
	\begin{split}
		\operatorname{HSC}^{0}(v)
		-\frac1{25}\left(7\mid\!\tau\!\mid^2+
		8\,\frac{\mid\!\tau(v)\!\mid^2}{\mid\!v\!\mid^2}\right)
		& \geq \kappa+
		\left(t^2+\frac8{25}\right)
		\left(\mid\!\tau\!\mid^2-\frac{\mid\!\tau(v)\!\mid^2}{\mid\!v\!\mid^2}\right)\\
		& \geq \kappa.
	\end{split}
	\]
	Thus \eqref{eq:directional-curvature-bound} holds and Theorem \ref{thm:directional-spectral-gap} applies. If one assumes only $\operatorname{HSC}^{t}\geq\kappa$, then \eqref{eq:main-pinching} holds with $\kappa-\frac35\mid\!\tau\!\mid_\infty^2$ in place of $\kappa$, which gives the second assertion.
\end{proof}

\begin{proof}[Proof of Corollary \ref{cor:closed-threshold}]
	If $f$ were nonconstant, the bound on $\psi'$ would contradict the strict inequality in Theorem \ref{thm:directional-spectral-gap}. Under \eqref{eq:main-pinching}, the same argument applies by Corollary \ref{cor:isotropic}.
\end{proof}

\section{Consequences of the curvature criteria} \label{sec:nonclosed-curvature-examples}

\subsection{}
In this part, we first construct (non-K\"{a}hler) Gauduchon metrics satisfying both criteria in a prescribed Aeppli class (e.g., \cite{DA}). We then perturb such a metric so that its $\operatorname{HSC}^{0}$ is arbitrarily negative on a prescribed line, while $\sigma_\omega$ remains positive.

\begin{proposition}\label{prop:nonclosed-positive-aeppli}
	Let $\left(X,\omega_0\right)$ be a compact connected K\"ahler surface with positive holomorphic sectional curvature, and set
	\[
	\kappa_0=\min_{x\in X,\, v\in T_x^{1,0}X,\, \mid v\mid=1}
	\HSC^{0}_{\omega_0}(v).
	\]
	For every $\kappa\in(0,\kappa_0)$, there are (non-K\"{a}hler) Gauduchon metrics $\omega_\varepsilon$ in the Aeppli class of $\omega_0$, converging smoothly to $\omega_0$ as $\varepsilon\to0$, such that
	\[
	\HSC^{0}_{\omega_\varepsilon}
	\geq\kappa+\frac35\mid\!\tau_\varepsilon\!\mid_{\omega_\varepsilon}^2.
	\]
	In particular, these metrics satisfy \eqref{eq:directional-curvature-bound} and $\sigma_{\omega_\varepsilon}\geq\kappa$.
\end{proposition}

\begin{proof}
	In a holomorphic coordinate ball choose a cutoff function $\chi$ supported in the ball and equal to one near its center. Put
	\[
	\alpha=\chi\mid\!z^1\!\mid^2\,\mathrm{d}z^2,\quad
	\zeta=\bar\partial\alpha+\partial\bar\alpha,
	\]
	and extend the forms by zero. The real (1,1)-form $\zeta$ is $\partial\bar\partial$-closed and Aeppli exact. Near the center,
	\[
	\mathrm{d}\zeta=\mathrm{d}z^1\wedge\mathrm{d}\bar z^1\wedge\mathrm{d}z^2+
	\mathrm{d}\bar z^1\wedge\mathrm{d}z^1\wedge\mathrm{d}\bar z^2\ne0.
	\]
	Hence, for sufficiently small $\varepsilon>0$,
	\[
	\omega_\varepsilon=\omega_0+\varepsilon\zeta
	\]
	is a (non-K\"{a}hler) Gauduchon metric, and its Aeppli class is unchanged.
	
	Smooth convergence of $\omega_\varepsilon$ gives uniform convergence of $\HSC^{0}$ on the $\omega_0$-unit tangent sphere bundle. Since $\omega_0$ is K\"ahler, $\tau_\varepsilon\to0$ uniformly. The asserted inequality follows from $\kappa<\kappa_0$.
	It implies \eqref{eq:directional-curvature-bound} since
	\[
	\frac{\mid\!\tau(v)\!\mid^{2}}{\mid\!v\!\mid^{2}}\leq\mid\!\tau\!\mid^{2}.
	\]
	Moreover, \eqref{eq:effective-quartic-potential} gives
	\[
	\mathcal H_{\omega_\varepsilon}
	\geq\kappa+\frac1{10}\mid\!\tau_\varepsilon\!\mid_{\omega_\varepsilon}^2
	\geq\kappa,
	\]
	so the variational characterization \eqref{eq:curvature-form-bottom} gives $\sigma_{\omega_\varepsilon}\geq\kappa$.
\end{proof}

For example, on $\mathbb{CP}^2$ with the Fubini-Study metric normalized by $\operatorname{HSC}=2$, one may take $\kappa=1$. The proposition also applies in every K\"ahler class on each Hirzebruch surface $\mathbb P\left(\mathcal O\oplus\mathcal O(a)\right)\to\mathbb{CP}^1$: for $a\geq1$, this follows from Yang-Zheng \cite{YangZheng}, and for $a=0$ from the product of suitably scaled Fubini-Study metrics.

\begin{proposition}\label{prop:nonclosed-local-negative}
	Let $\omega_0$ be a Gauduchon metric on a compact connected complex surface $X$. Fix
	\[
	p\in X,\quad v\in(T_p^{1,0}X\setminus\{0\}),\quad B>0. 
	\]
	There are Gauduchon metrics $\omega_\varepsilon$ in the Aeppli class of $\omega_0$ such that $\omega_\varepsilon\rightarrow\omega_0$ in $C^1$, and
	\begin{equation}\label{11074}
		\begin{split}
			\mathrm{d}\omega_\varepsilon &=\mathrm{d}\omega_0,\quad
			\Int\omega_\varepsilon^{2}=\Int\omega_0^{2}, \\
			\HSC^{t}_{\omega_\varepsilon}(v)&<-B, \quad\text{for every } t\in\mathbb R,\\
			\sigma_{\omega_\varepsilon}&=\sigma_{\omega_0}+O(\varepsilon^2).
		\end{split}
	\end{equation}
	Since $\mathrm{d}\omega_\varepsilon=\mathrm{d}\omega_0$, the perturbed metrics are non-K\"ahler whenever $\omega_0$ is non-K\"ahler.
	
	Moreover, if $\sigma_{\omega_0}>0$, then for sufficiently small $\varepsilon>0$,
	\[
	\lambda_1\left(\omega_\varepsilon\right)
	>\sigma_{\omega_\varepsilon}\geq\frac12\sigma_{\omega_0}>0,
	\]
	although \eqref{eq:directional-curvature-bound} fails for every $\kappa>0$.
\end{proposition}

\begin{proof}
	Choose holomorphic coordinates centered at $p$ in which
	\[
	g_{0,\,i\bar j}(p)=\delta_{ij},\quad
	\mathbb Cv=\mathrm{span}\left\{\frac{\partial}{\partial z^1}\right\}.
	\]
	Take $\phi\in C^\infty_c(\mathbb C^2,\mathbb R)$ supported in the unit ball and equal to one near the origin. For
	fixed $A>0$ and small $\varepsilon>0$, set
	\[
	u_\varepsilon(z)=A\mid\!z^1\!\mid^4\, \phi\left(\frac{z}{\varepsilon}\right),\quad
	\omega_\varepsilon=\omega_0+\sqrt{-1}\,\partial\bar\partial u_\varepsilon.
	\]
	Extend $u_\varepsilon$ by zero outside the chart. For $0\leq k\leq4$, its coordinate derivatives satisfy
	\[
	\mid\!D^k u_\varepsilon\!\mid\leq C_A\varepsilon^{4-k}.
	\]
	Thus $\omega_\varepsilon$ is positive for small $\varepsilon$, converges in $C^1$ to $\omega_0$, and has uniformly bounded Chern curvature. Direct calculation yields the first two identities in \eqref{11074}.
	
	At $p$ the metric and its first derivatives are unchanged. By \eqref{09221},
	\[
	\left(R_{\omega_\varepsilon}\right)_{1\bar1 1\bar1}(p)
	=\left(R_{\omega_0}\right)_{1\bar1 1\bar1}(p)-4A.
	\]
	The first coordinate vector is still unit at $p$. Choose $A>0$ so that $\HSC^{0}_{\omega_\varepsilon}(v)<-B$. By \eqref{eq:rankone},
	\[
	\HSC^{t}_{\omega_\varepsilon}(v)\leq\HSC^{0}_{\omega_\varepsilon}(v)<-B
	\]
	for every $t\in\mathbb R$. Thus the same perturbation works for all $t$. We keep $A$ fixed as $\varepsilon\to0$. In particular, the constant in the $O(\varepsilon^2)$ estimate in \eqref{11074} may depend on $B$.
	
	It remains to compare the Rayleigh quotients in \eqref{eq:curvature-form-bottom}. We use $\omega_0$ as a background metric in the following estimates. The potentials $\mathcal H_{\omega_\varepsilon}$ and
	$\mathcal H_{\omega_0}$ are uniformly bounded and agree outside the coordinate ball of radius $\varepsilon$. Since $X$ has real dimension 4,
	\[
	\left\lvert\mathcal H_{\omega_\varepsilon}-\mathcal H_{\omega_0}\right\rvert_{L^2}
	\leq C\varepsilon^2.
	\]
	The metric and its relative volume density differ from their background values by $O(\varepsilon^2)$ in $C^0$. Consequently, H\"older's inequality and the Sobolev inequality $W^{1,2}\hookrightarrow L^4$ give, for every real $u\in H^1$,
	\[
	\left\lvert\mathcal E_{\omega_\varepsilon}(u)-\mathcal E_{\omega_0}(u)\right\rvert
	\leq C\varepsilon^2
	\left(\mid\!\mathrm{d}u\!\mid_{L^2}^2+\mid\!u\!\mid_{L^2}^2\right).
	\]
	Indeed, the change in the potential at fixed volume contributes at most
	\[
	\left\lvert\mathcal H_{\omega_\varepsilon}-\mathcal H_{\omega_0}\right\rvert_{L^2}
	\mid\!u\!\mid_{L^4}^2.
	\]
	The other terms follow from metric and volume comparison.
	
	Take
	\[
	L=1+\mid\!\mathcal H_{\omega_0}\!\mid_\infty.
	\]
	The shifted form
	\[
	\mathcal E_{\omega_0}(u)+L\int_X u^2
	\]
	controls $\mid\!u\!\mid_{H^1}^2$. Hence the corresponding shifted form for $\omega_\varepsilon$ is bounded above and below by $1+C\varepsilon^2$ and $1-C\varepsilon^2$ times the shifted form for $\omega_0$. The same comparison holds for the two squared $L^2$ norms. Taking infima of the shifted Rayleigh quotients yields
	\[
	\sigma_{\omega_\varepsilon}+L=\sigma_{\omega_0}+L+O(\varepsilon^2).
	\]
	If $\sigma_{\omega_0}>0$, the spectral bound now follows from Theorem \ref{11062}. The negative curvature on $\mathbb Cv$ excludes \eqref{eq:directional-curvature-bound} for every positive $\kappa$.
\end{proof}

Taking an arbitray (non-K\"{a}hler) Gauduchon metric in Proposition \ref{prop:nonclosed-positive-aeppli} as a background metric, the second construction therefore separates the two criteria on the same complex surface and in the same Aeppli class: the integral criterion gives a positive spectral bound even when the pointwise curvature hypothesis of Theorem \ref{thm:directional-spectral-gap} fails.

\subsection{}
The preceding examples concern non-K\"{a}hler metrics on complex surfaces that admit K\"ahler metrics. We next consider a family of surfaces that admit no K\"ahler metric and for which the local curvature and torsion data remain fixed while the first eigenvalue tends to zero.

Let $\left(X, g\right)$ be a compact connected Riemann surface with area form $\omega_g$ such that $[\frac{\omega_g}{\pi}]$ is an integral class. Choose a principal circle bundle $\pi:P\to X$ with connection form $\eta$ satisfying
\[
\eta\left(\xi\right)=1,\quad \mathrm{d}\eta=2\,\pi^*\omega_g,
\]
where $\xi$ generates the circle action with period $2\pi$. For $\ell>0$, set
\[
X_\ell=\left(\mathbb R/\ell\mathbb Z\right)\times P, \quad
g_\ell=2\left(\mathrm{d}y^2+\eta^2+\pi^*g\right).
\]
Equip $X_\ell$ with the complex structure for which
\[
\alpha=\mathrm{d}y+\sqrt{-1}\,\eta
\]
and the pullbacks of local (1,0)-forms on $X$ are of type (1,0). This structure is integrable, since $\mathrm{d}\alpha=2\sqrt{-1}\,\pi^*\omega_g$ has type $(1,1)$. The fundamental form satisfies
\[
\omega_\ell=2\,\mathrm{d}y\wedge\eta+2\,\pi^*\omega_g, \quad
\mathrm{d}\omega_\ell=-2\,\mathrm{d}y\wedge\omega_\ell.
\]
Thus the Lee form is $-2\,\mathrm{d}y$. It is parallel for the product metric, so $\omega_\ell$ is Vaisman and hence Gauduchon. Since the circle bundle has nonzero Euler class, the Gysin sequence gives
\[
b_1\left(X_\ell\right)=1+2g(X).
\]
In particular, $X_\ell$ admits no K\"ahler metric.

The Chern holomorphic sectional curvature vanishes on the complex Lee line. Indeed, $\partial_y$ and $\xi=J\partial_y$ generate holomorphic isometries. In a local holomorphic coordinate $w$ for which
\[
\frac{\partial}{\partial w}=\frac12\left(\partial_y-\sqrt{-1}\,\xi\right),
\]
the metric coefficients are independent of both $\Real w$ and $\Imag w$. By \eqref{09221}, $R_{w\bar w w\bar w}=0$. Thus \eqref{eq:directional-curvature-bound} fails for every $\kappa>0$. The integral criterion gives no positive lower bound here either: the torsion trace is $\tau=-\alpha$, with $\mid\!\tau\!\mid^2=1$. By \eqref{eq:effective-quartic-potential} and \eqref{eq:curvature-form-bottom},
\[
\mathcal H_{\omega_\ell}\leq-\frac{1}{2}, \quad \sigma_{\omega_\ell}\leq-\frac{1}{2}.
\]

Finally, all the $X_\ell$ have the same Hermitian covering space $\mathbb R\times P$. Varying $\ell$ therefore leaves the local curvature and torsion unchanged. On the other hand, the function
\[
u_\ell(y)=\cos\left(\frac{2\pi}{\ell}y\right)
\]
descends to $X_\ell$ and satisfies
\[
\Delta_{\mathrm d}u_\ell=\frac{2\pi^2}{\ell^2}u_\ell.
\]
Consequently, as $\ell\to\infty$,
\[
\lambda_1\left(g_\ell\right)\leq\frac{2\pi^2}{\ell^2}.
\]
In this degenerate curvature regime, the Lee period is a global parameter that the local Hermitian geometry does not detect.

\begin{corollary}\label{cor:gs-positive-criteria-rationality}
	Let $(X, \omega)$ be a compact connected Gauduchon surface. If $\omega$ satisfies \eqref{eq:directional-curvature-bound} for some $\kappa>0$, or if $\sigma_\omega>0$, then $X$ is rational. Conversely, every rational surface admits a (non-K\"{a}hler) Gauduchon metric satisfying both conditions.
\end{corollary}

\begin{proof}
	Under either hypothesis, the first positive eigenvalues of all connected finite unramified covers of $X$, equipped with the pullback metric, have a common positive lower bound. Indeed, the directional curvature inequality pulls back unchanged, so Theorem \ref{thm:directional-spectral-gap} gives the same bound on every cover. For the integral criterion, the potential $\mathcal H_\omega$ also pulls back. A positive ground state of
	\[
	\frac14\Delta_{\mathrm d}+\mathcal H_\omega
	\]
	therefore lifts to a positive eigenfunction on the cover. On a compact connected cover such a eigenfunction has the lowest eigenvalue, so the spectral bottom is still $\sigma_\omega$. Now we can apply Theorem \ref{11062}.
	
	Morepver, we have $b_1(X)=0$. Otherwise, choose a smooth real closed one-form $\alpha$ whose periods form $\mathbb Z$. Reducing its period homomorphism modulo $m$ gives a connected cyclic cover
	\[
	\pi_m:X_m\to X
	\]
	of degree $m$. Its pullback has a primitive
	\[
	F_m:X_m\to\mathbb R/m\mathbb Z
	\]
	modulo $m$, with $\mathrm{d}F_m=\pi_m^*\alpha$. For $m\geq2$, a generator of the deck group can be chosen to translate $F_m$ by one. Thus
	\[
	u_m=\exp\left(2\pi\sqrt{-1}\, \frac{F_m}{m}\right)
	\]
	has integral zero. Applying the Rayleigh principle to its real and imaginary parts yields
	\[
	\lambda_1\left(X_m\right)
	\leq
	\frac{\int_{X_m}\mid\!\mathrm du_m\!\mid^2}{\int_{X_m}\mid\!u_m\!\mid^2}
	=
	\frac{4\pi^2}{m^2}\frac{\int_X\mid\!\alpha\!\mid^2\omega^{2}}{\int_X\omega^{2}}.
	\]
	Letting $m\to\infty$ contradicts the uniform lower bound.
	
	It remains to determine the Kodaira dimension. By \eqref{10111} and \eqref{11022},
	\[
	s-\hat{s}=-g^{k\bar l}\nabla_{\bar l}T_k=\mid\!\tau\!\mid^2,
	\]
	The average of $\operatorname{HSC}^{0}(v)$ over the unit sphere in $T^{1,0}_xX$, with normalized invariant measure, is $\frac{1}{6}\left(s+\hat{s}\right)$. The average of $v^i\bar v^jv^k\bar v^l$ is (cf. \cite{XKY})
	\[
	\frac{1}{6}\left(g^{i\bar{j}}g^{k\bar{l}}+g^{i\bar{l}}g^{k\bar{j}}\right).
	\]
	Under \eqref{eq:directional-curvature-bound} with $\kappa>0$, this average is positive, whence $s>0$. For the other hypothesis, the minimum of the holomorphic sectional curvature is bounded above by its spherical average. It follows that
	\[
	\mathcal H_\omega\leq\frac{1}{6}\left(s+\hat s\right)-\frac{1}{2}\mid\!\tau\!\mid^2
	=\frac{1}{3}\left(s-2\mid\!\tau\!\mid^2\right).
	\]
	Testing the Rayleigh quotient defining $\sigma_\omega$ with the constant function gives
	\[
	\sigma_\omega
	\leq\frac{\int_X\mathcal H_\omega\, \omega^{2}}{\int_X\omega^{2}}
	\leq\frac{\int_X\left(s-2\mid\!\tau\!\mid^2\right)\omega^{2}}{3\int_X\omega^{2}}.
	\]
	Hence $\int_Xs>0$ in either case. Thus the Gauduchon degree of the canonical bundle is
	\[
	\deg_\omega K_X=-\frac{1}{2\pi}\int_Xs<0.
	\]
	
	If a holomorphic line bundle $L$ admits a nonzero section with zero divisor $D$, the Poincar\'{e}–Lelong formula gives
	\[
	\deg_\omega L=\int_D\omega\ge0.
	\]
	Thus no positive power of $K_X$ has a nonzero holomorphic section, and $\kappa(X)=-\infty$. By the well-known Enriques-Kodaira classification, and taking $b_1$ into account, $X$ is rational.
	
	Conversely, if $X$ is rational, Zhang's Theorem 1.4 in \cite{ShiZhang} gives a K\"ahler metric of positive holomorphic sectional curvature on $X$. Applying Proposition \ref{prop:nonclosed-positive-aeppli} to this metric gives a (non-K\"{a}hler) Gauduchon metric satisfying \eqref{eq:directional-curvature-bound} with $\kappa>0$, $\sigma_\omega>0$.
\end{proof}

On every rational surface, the converse part of Corollary \ref{cor:gs-positive-criteria-rationality} provides a (non-K\"{a}hler) Gauduchon metric satisfying both positive criteria. Applying Proposition \ref{prop:nonclosed-local-negative} to this metric gives another metric in the same Aeppli class, with the same exterior derivative and total volume, for which $\sigma_\omega>0$ but \eqref{eq:directional-curvature-bound} fails for every $\kappa>0$. Thus the two criteria distinguish metrics even on a fixed rational surface.

\begin{remark}
	There is no uniform positive lower bound for $\lambda_1$ depending only on a positive lower bound for $\HSC^{0}$ of a compact Gauduchon surface. Indeed, Lin and Ren \cite{LinRenHopf} construct a homogeneous Hermitian metric $\omega$ with positive Chern real bisectional curvature on a Hopf surface $X$ with $\pi_1(X)\simeq\mathbb Z$. Homogeneity gives
	\[
	\sqrt{-1}\,\partial\bar\partial\omega=c\,\omega^2
	\]
	for a constant $c$. Integration gives $c=0$, so $\omega$ is Gauduchon. Taking a single
	nonzero weight in the definition of real bisectional curvature gives positive $\HSC^{0}$. By compactness, $\HSC^0_\omega\ge h$ for some $h>0$.
	
	For the connected cyclic cover
	\[
	\pi_m:\, X_m\to X
	\]
	of degree $m$, the pullback metric satisfies the same curvature bound. If $\alpha$ is a closed real one-form on $X$ with period group $\mathbb Z$, the estimate in the proof of Corollary \ref{cor:gs-positive-criteria-rationality} gives, for $m\ge2$,
	\[
	\lambda_1\left(X_m,\pi_m^*\omega\right)
	\le \frac{4\pi^2}{m^2}
	\frac{\displaystyle\int_X\mid\!\alpha\!\mid_\omega^2\,\omega^2}
	{\displaystyle\int_X\omega^2}\rightarrow 0.
	\]
	Thus positive $\HSC^{0}$ alone does not give a uniform spectral gap on Gauduchon surfaces.
\end{remark}

\end{document}